\documentclass[a4paper]{amsart}
\usepackage{amssymb}
\usepackage{verbatim}
\theoremstyle{plain}
  \newtheorem{thm}{Theorem}[section]

  \newtheorem{conj}[thm]{Conjecture}
  
  \newtheorem{obs}[thm]{Observation}
  \newtheorem*{obs*}{Observation}
\theoremstyle{definition}

\theoremstyle{remark}
  \newtheorem{rem}[thm]{Remark}
  
  \newtheorem*{ack}{Acknowledgments}

\newcommand{\C}{\mathbb{C}}

\newcommand{\Q}{\mathbb{Q}}

\newcommand{\Li}{\operatorname{Li}}
\newcommand{\arccosh}{\operatorname{arccosh}}

\renewcommand{\Re}{\operatorname{Re}}
\renewcommand{\Im}{\operatorname{Im}}
\numberwithin{equation}{section}
\allowdisplaybreaks
\begin{document}
\title
{A version of the volume conjecture}
\author{Hitoshi Murakami}
\address{
Department of Mathematics,
Tokyo Institute of Technology,
Oh-okayama, Meguro, Tokyo 152-8551, Japan
}
\email{starshea@tky3.3web.ne.jp}
\date{\today}
\begin{abstract}
We propose a version of the volume conjecture that would relate a certain limit of the colored Jones polynomials of a knot to the volume function defined by a representation of the fundamental group of the knot complement to the special linear group of degree two over complex numbers.
We also confirm the conjecture for the figure-eight knot and torus knots.
This version is different from S.~Gukov's because of a choice of polarization.
\end{abstract}
\thanks{This research is partially supported by Grant-in-Aid for Scientific
Research (B) (15340019).}
\maketitle
\section{Introduction}
For a knot $K$ in the three-sphere $S^3$, one can define the colored Jones polynomial $J_N(K;t)$ as the quantum invariant corresponding to the $N$-dimensional irreducible representation of the Lie algebra $sl(2;\C)$
\cite{Jones:BULAM385,Kirillov/Reshetikhin:1989}
(see also \cite{Kirby/Melvin:INVEM1991}).
\par
The volume conjecture \cite{Kashaev:LETMP97,Murakami/Murakami:ACTAM12001} states that the limit of
$\log\left(J_N\bigl(K;\exp(2\pi\sqrt{-1})\bigr)\right)/N$ would determine the simplicial volume of the knot complement $S^3\setminus{K}$.
In \cite{Murakami/Yokota:JREIA} Y.~Yokota and the author proved that for the figure-eight knot $E$ and a complex number $r$ the limit
$\log\left(J_N\bigl(E;\exp(2\pi{r}\sqrt{-1})\bigr)\right)/N$ also exists and defines the volume for the three-manifold obtained from $S^3$ by certain Dehn surgery if $r$ is close to $1$.
\par
In this paper we will show that a similar phenomenon appears for torus knots, which are not hyperbolic.
We also propose a version of the volume conjecture which relates the limit of
$\log\left(J_N\bigl(K;\exp(2\pi{r}\sqrt{-1})\bigr)\right)/N$ to the volume function corresponding to a representation of $\pi_1(S^3\setminus{K})$ to
$SL(2;\C)$ such that ratio of the eigenvalues of its image of the meridian is
$\exp\left(2\pi{r}\sqrt{-1}\right)$.
Our version of generalization of the volume conjecture is different from S.~Gukov's \cite{Gukov:COMMP2005} due to a choice of polarization.
\begin{ack}
The author would like to thank S.~Kojima and T.~Yoshida for helpful conversations.
Thanks are also due to S.~Garoufalidis, who pointed out
Remarks~\ref{rem:Garoufalidis1} and \ref{rem:Garoufalidis2}, and to S.~Gukov,
who described the difference between his conjecture and the current one.
\end{ack}
%%%%%%%%%%%%%%%%%%%%%%%%%%%%%%%%%%%%%%%%%%%%%%%%%%%%%%%%%%%%%%
\section{A conjecture}
Let $K$ be a knot in $S^3$.
Denote by $J_N(K;q)$ the colored Jones polynomial associated to the $N$-dimensional irreducible representation of the Lie algebra $sl(2,\C)$
\cite{Kirillov/Reshetikhin:1989}.
We normalize it so that $J_2(K;q)$ is the Jones polynomial \cite{Jones:BULAM385} and that $J_N(U;q)=1$, where $U$ is the unknot.
\par
We propose the following conjecture as a version of the volume conjecture.
\begin{conj}[Parameterized Volume Conjecture]\label{conj:parameterization}
There exists an open subset $\mathcal{O}_{K}$ of $\C$ such that for any
$u\in\mathcal{O}_{K}$ the following limit exists:
\begin{equation}\label{eq:limit}
  \lim_{N\to\infty}
  \frac{\log
        J_N\left(
             K;\exp\left(
                     \left(u+2\pi\sqrt{-1}\right)/N
                   \right)
           \right)}
       {N}.
\end{equation}
Moreover the function of $u$
\begin{equation}\label{eq:H}
  H(K;u)
  :=
  \left(u+2\pi\sqrt{-1}\right)
  \lim_{N\to\infty}
  \frac{\log
        J_N\left(
             E;\exp\left(
                     \left(u+2\pi\sqrt{-1}\right)/N
                   \right)
           \right)}
       {N}
\end{equation}
is analytic on $\mathcal{O}_{K}$.
If we put
\begin{equation*}
  v_{K}(u)
  :=
  2\frac{d\,H(K;u)}{d\,u}-2\pi\sqrt{-1},
\end{equation*}
then the following formula holds:
\begin{equation}\label{eq:volume}
  V(K;u)
  =
  \Im\bigl(H(K;u)\bigr)-\pi\Re(u)-\frac{1}{2}\Re(u)\Im\bigl(v_{K}(u)\bigr),
\end{equation}
where $V(K;u)$ is the volume function corresponding to the representation from $\pi_1(S^3\setminus{K})$ to $SL(2;\C)$ sending the meridian and the longitude to the elements the ratios of whose eigenvalues are $\exp(u)$ and $\exp\bigl(v_{K}(u)\bigr)$ respectively \cite[\S~4.5]{Cooper/Culler/Gillet/Long/Shalen:INVEM1994}.
\end{conj}
\begin{rem}
If $u$ is parameterized by a real number $t$, then $V\bigl(K;u(t)\bigr)$ satisfies the following differential equation from Schl{\"a}fli's formula:
\begin{multline}\label{eq:Schlafli}
  \frac{d\,V\bigl(K;u(t)\bigr)}{d\,t}
  \\
  =
  -\frac{1}{2}
  \left(
    \Re\bigl(u(t)\bigr)\frac{d\,\Im\left(v_{K}\bigl(u(t)\bigr)\right)}{d\,t}
   -\Re\left(v_{K}\bigl(u(t)\bigr)\right)\frac{d\,\Im\bigl(u(t)\bigr)}{d\,t}
  \right).
\end{multline}
See \cite[\S~5]{Neumann/Zagier:TOPOL85}
and \cite[\S~4.5]{Cooper/Culler/Gillet/Long/Shalen:INVEM1994}.
Note that we use the same convention for the meridian/longitude pair as in
\cite{Neumann/Zagier:TOPOL85}, which is different from that in
\cite{Cooper/Culler/Gillet/Long/Shalen:INVEM1994} and \cite{Gukov:COMMP2005}.
Note also that the right hand side of the equation in the last line of Page~62
of \cite{Cooper/Culler/Gillet/Long/Shalen:INVEM1994} should be multiplied by
four (see \cite[(5.6)]{{Gukov:COMMP2005}}).
\end{rem}
\begin{rem}\label{rem:Dehn}
If there exist coprime integers $p$ and $q$ satisfying
$pu+qv_{K}(u)=2\pi\sqrt{-1}$, then $u$ would define the $(p,q)$-Dehn surgery along $K$ \cite{Thurston:GT3M}.
If this is a hyperbolic manifold, $V(K;u)$ is its hyperbolic volume.
\end{rem}
\begin{rem}
The open set $\mathcal{O}_{K}$ may not contain $0$.
Therefore Conjecture~\ref{conj:parameterization} is {\em not} a generalization of the volume conjecture.
(Recall that the volume conjecture
\cite{Kashaev:LETMP97,Murakami/Murakami:ACTAM12001}
states that when $u=0$, then the limit
\eqref{eq:limit} gives the simplicial volume of the knot complement.)
In fact the case where $K$ is a torus knot, the limit \eqref{eq:limit} is not continuous at $0$ \cite{Kashaev/Tirkkonen:ZAPNS2000,Murakami:INTJM62004}.
See also \cite[Proposition~B.2]{Garoufalidis/Le:Kauffman}.
\end{rem}
\begin{rem}\label{rem:Garoufalidis1}
S.~Garoufalidis and T.~Le proved that if $u$ is close enough to
$-2\pi\sqrt{-1}$, then $J_N\left(K;\exp\bigl((u+2\pi\sqrt{-1})/N\bigr)\right)$
converges to $1/\Delta(K;\exp(u+2\pi\sqrt{-1})\bigr)$, where $\Delta(K;t)$
is the Alexander polynomial of $K$ \cite{Garoufalidis/Le:aMMR}.
(See also \cite{Murakami:2005} for the figure-eight knot.)
In this case the right hand side of \eqref{eq:limit} vanishes and so we have
$H(K;u)=0$ and $v_K(u)=-2\pi\sqrt{-1}$.
Therefore from \eqref{eq:volume} the volume function $V(K;u)$ vanishes.
This corresponds to the case where $u$ defines an abelian representation, whose
volume function is zero.
See \cite[Appendix B]{Garoufalidis/Le:Kauffman}.
So we exclude such a case in Conjecture~\ref{conj:parameterization}.
\end{rem}
\par
Note that the conjecture above is proved for the figure-eight knot by Yokota and the author \cite{Murakami/Yokota:JREIA}.
In fact we proved \cite[Corollary 2.4]{Murakami/Yokota:JREIA}
that for the figure-eight knot $E$
\begin{equation*}
  V(E;u)
  =
  \Im\bigl(H(E;u)\bigr)-\pi\Re(u)-\frac{1}{4}\Im\bigl(u\,v_E(u)\bigr)
  -\frac{\pi}{2}\operatorname{length}(\gamma),
\end{equation*}
where $\operatorname{length}(\gamma)$ is the length of the geodisic $\gamma$ added to complete the incomplete hyperbolic structure of $S^3\setminus{E}$ corresponding to $u$.
But since
$\operatorname{length}(\gamma)=-\Im\bigl(u\,\overline{v_E(u)}\bigr)/(2\pi)$
from \cite[(34)]{Neumann/Zagier:TOPOL85}, we have
\begin{equation*}
\begin{split}
  V(E;u)
  &=
  \Im\bigl(H(E;u)\bigr)-\pi\Re(u)-\frac{1}{4}\Im\bigl(u\,v_E(u)\bigr)
  +\frac{1}{4}\Im\bigl(u\,\overline{v_E(u)}\bigr)
  \\
  &=
  \Im\bigl(H(E;u)\bigr)-\pi\Re(u)-\frac{1}{2}\Re(u)\Im\bigl(v_E(u)\bigr)
\end{split}
\end{equation*}
as required.
\par
See also \cite[(4.1)]{Murakami:FUNDM2004}.
(The sign of $\pi\Re{u}$ in \cite[(4.1)]{Murakami:FUNDM2004} should be negative because the author used a wrong definition of the function $H$ in the old version of \cite{Murakami/Yokota:JREIA}.)
\par
Note also that Gukov uses a different {\em polarization}
in his generalization of the volume conjecture \cite[(5.12)]{Gukov:COMMP2005}.
It agrees with Conjecture~\ref{conj:parameterization} when $\Re(u)=0$.
The difference for $\Re(u)\ne0$ can be explained by a choice of polarization.
Details will be described in a forthcoming paper.
%%%%%%%%%%%%%%%%%%%%%%%%%%%%%%%%%%%%%%%%%%%%%%%%%%%%%%%%%%
\section{Proof}
In this section we prove Conjecture~\ref{conj:parameterization} for torus knots.
(To be honest, \eqref{eq:volume} is proved only up to a constant, that is, we prove \eqref{eq:Schlafli} instead.)
\par
Let $T(a,b)$ be the torus knot of type $(a,b)$, where $a$ and $b$ are coprime integers with $a>1$ and $b>1$.
Then the author proved in \cite[Theorem~1.1]{Murakami:INTJM62004} the following theorem.
\begin{thm}[\cite{Murakami:INTJM62004}]
Suppose that $|r|>1/(ab)$, $\Re{r}>0$ and $\Im{r}>0$, then
\begin{equation*}
  \lim_{N\to\infty}
  \frac{\log{J_N\left(T(a,b);\exp(2\pi r\sqrt{-1}/N)\right)}}{N}
  =
  \left(1-\frac{1}{2abr}-\frac{abr}{2}\right)\pi\sqrt{-1}.
\end{equation*}
\end{thm}
\par
Therefore the functions $H\bigl(T(p,q);u\bigr)$ and $v_{T(a,b)}(u)$ in Conjecture~\ref{conj:parameterization} are defined as follows:
\begin{align*}
  H\bigl(T(p,q);u\bigr)
  &=
  (u+2\pi\sqrt{-1})
  \left(
    1
    -\frac{1}{2ab\left(1+\frac{u}{2\pi\sqrt{-1}}\right)}
    -\frac{ab\left(1+\frac{u}{2\pi\sqrt{-1}}\right)}{2}
  \right)
  \pi\sqrt{-1}
  \\
  &=
  \frac{-(ab(u+2\pi\sqrt{-1})-2\pi\sqrt{-1})^2}{4ab}
  \\
  \intertext{and}
  v_{T(a,b)}(u)
  &=
  2\frac{d\;H\bigl(T(a,b);u\bigr)}{d\;u}-2\pi\sqrt{-1}
  =
  -ab(u+2\pi\sqrt{-1})
\end{align*}
for $u$ with $|u+2\pi\sqrt{-1}|>2\pi/(ab)$, $\Re(u)<0$ and $\Im(u)>-2\pi$.
\par
So the volume function $V\bigl(T(a,b);u\bigr)$ is
\begin{equation*}
\begin{split}
  V\bigl(T(a,b);u\bigr)
  &=
  \Im{\bigl(H(T(a,b);u\bigr)}-\pi\Re(u)-\frac{1}{2}\Re(u)\Im(v_{T(a,b)}(u))
  \\
  &=
  -\frac{1}{4}ab\Im(u^2)
  -ab\pi\Re(u)
  +\pi\Re(u)
  \\
  &\phantom{=}
  -\pi\Re(u)
  -\frac{1}{2}\Re(u)(-ab\Im(u)-2ab\pi)
  \\
  &=
  0.
\end{split}
\end{equation*}
\par
On the other hand the right hand side of \eqref{eq:Schlafli} equals
\begin{equation*}
  -\frac{1}{2}
  \left(
    \Re\bigl(u(t)\bigr)
    (-ab)
    \frac{d\,\Im\bigl(u(t)\bigr)}{d\,t}
   +ab\Re\bigl(u(t)\bigr)
    \frac{d\,\Im\bigl(u(t)\bigr)}{d\,t}
  \right)
  =
  0,
\end{equation*}
which proves \eqref{eq:Schlafli}.
\par
This confirms Conjecture~\ref{conj:parameterization} for the torus knot $T(a,b)$.
\begin{rem}
The volume function for a torus knot would be zero.
See \cite[Appendix~B]{Garoufalidis/Le:Kauffman}.
\end{rem}
\begin{rem}
Note that the pair $\left(u,v_{T(a,b)}(u)\right)$ satisfies the equality:
\begin{equation*}
  (-1)\times{u}+\left(\frac{-1}{ab}\right)\times{v_{T(a,b)}(u)}=2\pi\sqrt{-1}.
\end{equation*}
So the corresponding geometric object would be the generalized Dehn surgery along the torus knot $T(a,b)$ with invariant $\bigl(1,1/(ab)\bigr)$, or the $ab$-Dehn surgery with cone-angle $2ab\pi$
\cite[Chapter~4, \S~4.5]{Thurston:GT3M} (see also Remark~\ref{rem:Dehn}).
It would be interesting that the $ab$-Dehn surgery along the torus knot $T(a,b)$ is reducible \cite{Moser:PACJM1971}.
\end{rem}
%%%%%%%%%%%%%%%%%%%%%%%%%%%%%%%%%%%%%%%%%%%%%%%%%%%%%%%%%%%%%%%%%%%%%%%
\section{Comments}
In \cite[(5.29)]{Gukov:COMMP2005}, Gukov proposed the following conjecture.
\begin{conj}\cite{Gukov:COMMP2005}\label{conj:Gukov}
Let $K$ be a knot in the three-sphere.
For $a\in\C\setminus\Q$, define the function $l(a)$ as follows:
\begin{equation}\label{eq:lambda}
  l(a):=
  -\frac{d}{d\,a}
  \left\{
    a\lim_{N\to\infty}
    \frac{\log{J_N\left(K;\exp(2\pi{a}\sqrt{-1}/N)\right)}}{N}
  \right\}.
\end{equation}
Then the pair $\left(\exp\bigl(l(a)\bigr),-\exp(\pi{a}\sqrt{-1})\right)$ is a zero of the $A$-polynomial of $K$ introduced in
\cite{Cooper/Culler/Gillet/Long/Shalen:INVEM1994}.
\end{conj}
Using our parameterization $u:=2\pi{a}\sqrt{-1}-2\pi\sqrt{-1}$, we have
$l(a)=-v_{K}(u)/2-\pi\sqrt{-1}$.
Then Conjecture~\ref{conj:Gukov} states that the pair
$\left(-\exp\bigl(-v_{K}(u)/2\bigr),\exp(u/2)\right)$
is a zero of the $A$-polynomial.
\par
In the case of the figure-eight knot, we can prove this.
For the torus knot $T(a,b)$, either
$\left(-\exp\bigl(-v_{K}(u)/2\bigr),\exp(u/2)\right)$ or 
$\left(-\exp\bigl(v_{K}(u)/2\bigr),\exp(u/2)\right)$ is a zero of the $A$-polynomial.
(See, for example, \cite[Example~4.1]{XZhang:2004} for the $A$-polynomials of torus knots.)
This would depend on how we choose a meridian/longitude pair.
\begin{rem}\label{rem:Garoufalidis2}
If $|a|$ is small enough the right hand side of \eqref{eq:lambda} vanishes \cite{Garoufalidis/Le:Kauffman}.
This corresponds to the $(\mathfrak{l}-1)$-factor of the $A$-polynomial
\cite[2.5]{Cooper/Culler/Gillet/Long/Shalen:INVEM1994}.
See Remark~\ref{rem:Garoufalidis1}.
\end{rem}
The function $H(K;u)$ defined by \eqref{eq:H} is a kind of potential function introduced in \cite[Theorem~3]{Neumann/Zagier:TOPOL85}.
(To be more precise, $\Phi(u)=4H(K;u)-4\pi{u}\sqrt{-1}$, where $\Phi(u)$ is Neumann--Zagier's potential function.)
The following observation indicates a relation between $H(K;u)$ and the Alexander polynomial.
\begin{obs}
Let $\Delta(K;t)$ be the Alexander polynomial of a knot $K$.
For the figure-eight knot and torus knots, the equations $H(K;u)=0$ and $\Delta\bigl(K;\exp(u)\bigr)=0$ share a root.
Note that here we regard $H(K;u)$ as a function defined on the whole complex plane.
\end{obs}
\begin{proof}
First note that $\Delta(E;t)=-t^2+3t-1$ and
$\Delta\bigl(T(a,b);t\bigr)
=\left(t^{ab}-1\right)(t-1)/\left(t^a-1\right)\left(t^b-1\right)$.
\par
For the figure-eight knot $E$, we have from \cite{Murakami/Yokota:JREIA}
\begin{equation*}
  H(E;u)
  =
  \Li_2(y^{-1}\exp(-u))-\Li_2(y\exp(-u))+\left(\log(-y)+\pi\sqrt{-1}\right)u,
\end{equation*}
where $y$ is a solution to the equation $y+y^{-1}=2\cosh(u)-1$.
Put $u:=\pm\arccosh(3/2)$.
Then $y=1$ and so $H(E;u)$ vanishes.
Note that we take the branch of $\log$ so that $\log(-1)=-\pi\sqrt{-1}$.
It is easy to see that $\exp(\pm\arccosh(3/2))$ are the roots of $\Delta(E;\exp(u))=0$.
\par
For the torus knot $T(a,b)$ the zero of $H\bigl(T(a,b);u\bigr)$ is
$2\pi\sqrt{-1}/(ab)-2\pi\sqrt{-1}$.
It is also a zero of $\Delta(T(a,b);\exp(u))$.
\end{proof}
%%%%%%%%%%%%%%%%%%%%%%%%%%%%%%%%%%%%%%%%%%%%%%%%%%%%%%%%%%%%%%%%%%%%
\bibliography{mrabbrev,hitoshi}
\bibliographystyle{hamsplain}
\end{document}